\documentclass[12pt]{article}
\usepackage{amsmath}
\usepackage{latexsym}
\usepackage{amssymb}
\newtheorem{la}{Lemma}
\newtheorem{Defn}{Definition}
\newtheorem{Remark}{Remark}
\newtheorem{Note}{Note}

\newtheorem{Example}{Example}
\newtheorem{Examples}{Examples}
\newtheorem{Problems}{Problems}

\newtheorem{Problem}{Problem}
\newtheorem{Number}{\!\!}

\newenvironment{proof}{{\noindent\bf Proof.}}%
                  {\nopagebreak\hspace*{\fill}$\Box$\medskip\medskip\par}

\newcommand{\wb}{\overline}

\newcommand{\at}{\symbol{'100}}

\newcommand{\n}{\rm}

\newcommand{\mto}{\mapsto}

\newcommand{\N}{{\mathbb N}}

\newcommand{\Q}{{\mathbb Q}}

\newcommand{\Z}{{\mathbb Z}}

\newcommand{\Aut}{\mbox{\n Aut}}

\newcommand{\sub}{\subseteq}

\newcommand{\im}{\mbox{\n im}}

\newcommand{\cB}{{\cal B}}

\newcommand{\cT}{{\cal T}}
\newcommand{\cK}{{\cal K}}

\newcommand{\dsemi}{\mbox{$\times\!$\rule{.15mm}{2mm}}\,}
\begin{document}
\begin{center}
{{\Large\bf Contraction groups for tidy automorphisms\\[2mm]
of totally disconnected groups}}\\[6mm]
{\bf Helge Gl\"{o}ckner\footnote{The present research
was supported by DFG grant 447 AUS-113/22/0-1
and ARC grant LX 0349209. The author thanks
U. Baumgartner and G.\,A. Willis for discussions
concerning the tidying procedure.}
}\vspace{5mm}
\end{center}
\begin{abstract}
\hspace*{-7.2 mm}
In this note, we show that results of U. Baumgartner and G.\,A. Willis
concerning contraction groups
of automorphisms of metrizable totally disconnected, locally compact groups
remain valid also in the non-metrizable
case, if one restricts attention to
automorphisms for which small tidy subgroups
exist.\vspace{2mm}
\end{abstract}
{\footnotesize{\em Key words}: locally compact group, totally disconnected
group, scale function, scale, tidy subgroup,
automorphism, contraction group, opposite subgroup,
metrizability\\[2mm]
{\em Classification}:
22D05, 
22D45.}\\[5mm]
Given an automorphism~$\alpha$ of a locally compact group~$G$,
its contraction group $U_\alpha$
is the group of all $x\in G$ such that $\lim_{n\to\infty}\alpha^n(x)=1$.
Contraction groups arise in the study of semistable
convolution semigroups on
second countable
groups~\cite{Haz}, \cite{Sie}.
Contraction groups in $p$-adic Lie groups have been investigated
in \cite{DaS}, and recently general results
for metrizable totally disconnected groups~$G$ were
obtained, using
the concept of a tidy subgroup (see \cite{Wi1}, \cite{Wi2})
as a tool~\cite{BaW}. In particular, it was shown in~\cite{BaW}
that $U_\alpha$ is closed if and only if $G$ has small subgroups
tidy for~$\alpha$, and various reformulations and
interesting consequences of this property were obtained.
Let us say that
{\em $\alpha$ is tidy\/} if $G$ has small
subgroups tidy for~$\alpha$.
While every automorphism of a $p$-adic Lie group
is tidy~\cite[Rem.\,3.33\,(2)]{BaW},
already for automorphisms of Lie groups
over local fields of positive characteristic tidiness is
an important additional
regularity property,
which is not always satisfied~\cite{SC2}.
In this setting,
the tidiness property characterizes those
analytic automorphisms whose scale can be calculated on the Lie algebra level~\cite{SC2}.
Presuming tidiness of the automorphisms
involved, it is frequently
possible to extend
results for $p$-adic Lie groups
(as in \cite{SCA}, \cite{GaW} or \cite{Wan})
to the case of positive characteristic~\cite{SC2}.\\[3mm]
In this note, we describe consequences of the tidiness property.
Based on techniques from~\cite{BaW},
we establish two lemmas which 
allow us to transfer many of the results from~\cite{BaW} 
to non-metrizable groups,
in the special
case of tidy automorphisms. We also show that an automorphism~$\alpha$
is tidy if and only if its restriction $\alpha|_{M_\alpha}$
to the ``Levi factor''
$M_\alpha:=\{x\in G\!: \mbox{$\alpha^\Z(x)$ is relatively compact}\}$
is tidy.\\[3mm]
Throughout the following, $G$ is a totally disconnected, locally compact
topological group and $\alpha\in \Aut(G)$.
By \cite[Prop.\,3.4]{BaW},
$U_\alpha$ is a normal subgroup of
$P_\alpha:=\{x\in G\!: \mbox{$\alpha^\N(x)$ is relatively compact}\}$.
Following~\cite{BaW}, we set $U_0:=\wb{U_\alpha}\,\cap\,\wb{U_{\alpha^{-1}}}$.
Given a subset $H\sub G$ and $x\in G$,
we write $\lim_{n\to\infty} \alpha^n(x)=1 \!\mod H$
if, for every identity neighbourhood $V\sub G$,
there exists $N\in \N$ such that $\alpha^n(x)\in VH$
for all $n\geq N$.\,\footnote{The definition in \cite{BaW}
only requires $\alpha^n(x)\in VHV$, but it is the above property which
is really used in the proof of {\em loc.\,cit.}
Theorem~3.8 (see p.\,9, line~2).}
We write $U_{\alpha/H}:=\{x\in G\!:
\lim_{n\to\infty} \alpha^n(x)=1\! \mod H\}$.\\[3mm]
We recall that a compact, open subgroup $V\sub G$ is called {\em tidy for~$\alpha$\/}
if
\begin{itemize}
\item[{\bf T1}]
$V=V_+V_-$, where $V_\pm:=\bigcap_{n\in \N_0}\alpha^n(V)$; \,and
\item[{\bf T2}]
The subgroups $V_{++}:=\bigcup_{n\in \N_0}\alpha^n(V_+)$
and $V_{--}:=\bigcup_{n\in \N_0}\alpha^{-n}(V_-)$
are closed in~$G$.
\end{itemize}
Note that $\alpha(V_+)\supseteq V_+$ and $\alpha(V_-)\sub V_-$. Thus $V_-\sub P_\alpha$.
Tidy subgroups always exists, and the index $s_G(\alpha):=[\alpha(V_+):V_+]$
is independent of the choice of tidy subgroup~$V$ (see \cite{Wi1}).
We abbreviate $V_0:=V_+\cap V_-$.
If $V$ is tidy, then
$P_\alpha\cap V=V_-$ by \cite[Rem.\,3.1]{BaW},
whence $P_\alpha$ is closed.
As a consequence,
also $M_\alpha=P_\alpha\cap P_{\alpha^{-1}}$ is closed,
and $V_0=V\cap M_\alpha$ if~$V$ is tidy \cite[La.\,9]{Wi1}.
If $H\leq G$ is an $\alpha$-stable closed subgroup
(i.e., $\alpha(H)=H$)
and $V\sub G$ is tidy for~$\alpha$, then there exists
a subgroup~$W\sub G$ tidy for~$\alpha$ such that
$W\sub V$ and $W\cap H$ is tidy for~$\alpha|_H$ \cite[La.\,4.1]{Wi2}.
Hence $\alpha|_H$ is tidy if so is~$\alpha$.
We now adapt \cite[Thm.\,3.8]{BaW}
to our setting (cf.\ \cite{HaS} and \cite[Thm.\,3.1]{DaS},
where analogous conclusions are achieved
for real, resp., $p$-adic Lie groups and compact~$H$).
\begin{la}\label{la1}
Let $H$ be an $\alpha$-stable
closed subgroup of~$G$. If 
$\alpha|_{M_\alpha}$ is tidy,
then
\begin{equation}\label{eq38}
U_{\alpha/H}\; =\; U_\alpha H\, .
\end{equation}
\end{la}
\begin{proof}
The inclusion $U_\alpha H\sub U_{\alpha/H}$ is trivial.
For the converse inclusion, let $x\in U_{\alpha/H}$;
we have to find $h\in H$ such that $xh\in U_\alpha$.
Let $O$ be a compact, open
subgroup of~$G$ tidy for~$\alpha$. By \cite[La.\,3.10]{BaW},
there exists $h_0\in H$ such that $xh_0\in U_{\alpha/H\cap O_0}$.
Since $O_0=M_\alpha\cap O$, after replacing $H$ with $H\cap O_0$
and $x$ with $xh_0$
we may assume now without loss of generality
that $H\sub M_\alpha$ and $H$ is compact.
Let $\cB$ be the set of compact open subgroups of~$G$.
For each $V\in \cB$, \cite[La.\,3.10]{BaW}
yields $h_V\in H$ and $N_V\in \N$ such that $\alpha^n(xh_V)\in V$
for all $n\geq N_V$. Since~$H$ is compact, the net $(h_V)_{V\in \cB}$
has a cluster point~$h\in H$.
We claim that
$\alpha^n(xh)\to 1$ as $n\to\infty$,
i.e., $xh\in U_\alpha$.
Given an identity neighbourhood
$Y\sub G$,
there exists an identity neighbourhood $W\sub G$
such that $WW\sub Y$;
since $M_\alpha$ has small subgroups tidy for~$\alpha|_{M_\alpha}$,
after shrinking~$W$ we may assume that $W\cap M_\alpha$
is a compact open subgroup of~$M_\alpha$
tidy for $\alpha|_{M_\alpha}$.
Then $H\cap W=H\cap M_\alpha\cap W=H\cap W_0$
is $\alpha$-stable.
By definition of~$h$ as a cluster point, there exists
$V\in \cB$ with $V\sub W$ such that $h_V^{-1}h\in W$
and thus $h_V^{-1}h\in H\cap W$.
For all $n\geq N_V$, we now obtain
\[
\alpha^n(xh)=
\alpha^n(xh_V)\alpha^n(h_V^{-1}h)\, \in \, V\alpha^n(H\cap W)=V (H\cap W)\sub  WW\sub  Y\,.
\]
Hence indeed $\alpha^n(xh)\to 1$.
\end{proof}
As a consequence of Lemma~\ref{la1} just established,
Proposition~3.7--Corollary~3.30 from \cite{BaW}
remain valid also for non-metrizable~$G$,
if one assumes in addition that $\alpha|_{M_\alpha}$ be tidy.
In fact, inspection shows that metrizability is never used directly
in the proofs of these results,
but only the validity of Eqn.\,(\ref{eq38}) from above.
The proof of Lemma~3.31
in {\em loc.\,cit.}, however,
requires metrizability of~$G$ for a second
reason. The lemma exploits an algorithm
for the construction of tidy subgroups described in~\cite[\S2]{FLA}.
This ``tidying procedure''
involves a certain subgroup, $K$,
whose definition we presently recall.
If~$G$ is metrizable,
then $K=\wb{U_\alpha\cap P_{\alpha^{-1}}}$
(see \cite[p.\,4]{FLA}),
and this equality is used essentially in~\cite{BaW}.
Unfortunately, it is not known whether
equality persists for non-metrizable~$G$.
Therefore, to transfer
\cite[La.\,3.31 and Thm.\,3.32]{BaW}
to our setting, we first need to
discuss~$K$ by hand, assuming that~$\alpha|_{M_\alpha}$ is tidy.
Let us recall the definition from~\cite[p.\,4]{FLA}.\\[3mm]
{\bf Definition.}
Write
$\cK_O :=
\{x\in M_\alpha : \, \mbox{$\alpha^n(x)\in O$ for all sufficiently large $n$}\}$
and $K_O:=\wb{\cK_O}$,
for each compact open subgroup $O\sub G$.
We define $K:=\bigcap\{K_O:\,\mbox{$O$ a compact open subgroup of~$G$}\}$.\\[3mm]
Note that $\cK_O= O_{--}\cap M_\alpha =  (O\cap M_\alpha)_{--}$
here.
\begin{la}\label{la2}
If $\alpha|_{M_\alpha}$ is tidy, then $K=U_0=\{1\}$.
\end{la}
\begin{proof}
Let us show first that $K \sub \bigcap \{W_{--} \! : W\in \cT(\alpha|_{M_\alpha})\}$,
where $\cT(\alpha|_{M_\alpha})$ is the set of
subgroups $W\sub M_\alpha$ tidy for~$\alpha|_{M_\alpha}$.
To this end, let $W\in \cT(\alpha|_{M_\alpha})$.
Then there exists a compact open subgroup $O\sub G$ such that $O\cap M_\alpha\sub W$.
Hence $\cK_O=(O\cap M_\alpha)_{--}\sub W_{--}$
and thus $K\sub K_O\sub W_{--}$, using that $W_{--}$ is closed. The assertion follows.
Since $M_\alpha=M_{\alpha|_{M_\alpha}}$,
we have $W_-=W=W_+$ for each $W\in \cT(\alpha|_{M_\alpha})$
by \cite[La.\,3.19]{BaW},
whence $\alpha(W)=W$ and thus $W_{--}=W$.
Summing up, we have
$K\sub \bigcap\cT(\alpha|_{M_\alpha})$.
But $\bigcap\cT(\alpha|_{M_\alpha})
=\{1\}$, as $\alpha|_{M_\alpha}$ is tidy.
Thus $K=\{1\}$.
Using \cite[Cor.\,3.27]{BaW}
and
$\wb{U_\alpha}\cap \wb{U_{\alpha^{-1}}} \sub P_\alpha\cap P_{\alpha^{-1}}=M_\alpha$,
we also get
$\{1\} =\bigcap \cT(\alpha|_{M_\alpha})
=\wb{U_{\alpha|_{M_\alpha}}}\,
\cap \, \wb{U_{\alpha^{-1}|_{M_\alpha}}}
=\wb{U_\alpha} \, \cap \, \wb{U_{\alpha^{-1}}}=U_0$.
\end{proof}
Stimulated by the results in~\cite{BaW}
(notably Prop.\,3.7, Prop.\,3.21 and Thm.\,3.32)
and \cite[Thm.\,3.5\,(iii)]{Wan},
we now formulate various useful conclusions:\\[3mm]
{\bf Theorem.}
{\em Let $G$ be a totally disconnected, locally
compact group.
Then $\alpha\in \Aut(G)$ is tidy if and only if $\alpha|_{M_\alpha}$
is tidy.
In this case, we have:}
\begin{itemize}
\item[\rm (a)]
$U_0=\{1\}$.
\item[\rm (b)]
{\em $U_\alpha$ and $U_{\alpha^{-1}}$ are closed.}
\item[\rm (c)]
{\em $U_\alpha\cap M_\alpha=U_{\alpha^{-1}}\cap M_\alpha=\{1\}$.}
\item[\rm (d)]
{\em Every compact open subgroup of~$G$ satisfying {\bf T1}
is tidy.}
\item[\rm (e)]
{\em $P_\alpha=U_\alpha\dsemi M_\alpha$ and $P_{\alpha^{-1}}=U_\alpha\dsemi M_{\alpha^{-1}}$
as topological groups.}
\item[\rm (f)]
{\em $U_\alpha M_\alpha U_{\alpha^{-1}}$ is an $\alpha$-stable, open subset
of~$G$ which contains every subgroup tidy for~$\alpha$.
The product map $\mu\!: U_\alpha\times M_\alpha\times
U_{\alpha^{-1}}\to U_\alpha M_\alpha U_{\alpha^{-1}}$, $(x,y,z)\mto xyz$ is a homeomorphism.}
\item[\rm (g)]
{\em $s_G(\alpha^{-1})=s_{P_\alpha}(\alpha^{-1}|_{P_\alpha})
=s_{U_\alpha}(\alpha^{-1}|_{U_\alpha})
=s_{V_{--}}(\alpha^{-1}|_{V_{--}})
$ for each subgroup $V\sub G$ tidy for~$\alpha$, where
each of $P_\alpha$, $V_{--}$ and $U_\alpha$ are $\alpha$-invariant closed subgroups
of~$P_\alpha$.}
\item[\rm (h)]
{\em $s_H(\alpha^{-1}|_H)=\Delta_H(\alpha^{-1}|_H)$
for each $\alpha$-invariant closed subgroup~$H$ of~$P_\alpha$,
where $\Delta_H\!: \Aut(H)\to\Q^+$
is the modular function.
In particular, $s_G(\alpha^{-1})=\Delta_{U_\alpha}(\alpha^{-1}|_{U_\alpha})$.}
\item[\rm (i)]
{\em If $q\!: G\to Q$ is a quotient morphism
with $\alpha(\ker(q))=\ker(q)$ and $\wb{\alpha}$ the induced automorphism
of~$Q$, then $q(U_\alpha)=U_{\wb{\alpha}}$.}\vspace{2mm}
\end{itemize}
\begin{proof}
We already know that tidiness of~$\alpha$ entails tidiness
of $\alpha|_{M_\alpha}$. We assume that $\alpha|_{M_\alpha}$
is tidy now and show that that this entails all of (a)--(i) as well as tidiness
of~$\alpha$. We have already seen that (a) holds;
(g), (h) and (i) are the analogues \cite[Prop.\,3.7 and Prop.\,3.21]{BaW}
in our setting, which are remain valid by Lemma~\ref{la1} above.

(b) Using (a) and \cite[Cor.\,3.30]{BaW}, we obtain
$\wb{U_\alpha}=U_0U_\alpha=\{1\}U_\alpha=U_\alpha$.
Replacing $\alpha$ with $\alpha^{-1}$ (whose restriction to $M_\alpha=M_{\alpha^{-1}}$
is tidy since so is~$\alpha$), we see that also $U_{\alpha^{-1}}$ is closed.

(c) By (a), (b) and \cite[La.\,3.29]{BaW},
we have $U_\alpha\cap M_\alpha=
\wb{U_\alpha}\cap M_\alpha=U_0=\{1\}$.
Replacing $\alpha$ with
$\alpha^{-1}$, we also obtain the second formula.

(d) Let $V\sub G$ be a compact open subgroup
satisfying {\bf T1}. We apply the tidying procedure
from \cite[\S2]{FLA} to create a subgroup
tidy for~$\alpha$. Since~$V$ satisfies {\bf T1},
Step~1 of the algorithm outputs~$V$.
In Step~2a, one defines~$K$ as above;
thus $K=\{1\}$. In Step~3a,
one defines $V'':=\{v\in V\!: xvx^{-1}\in VK\}$;
then $V''=V$ as $K=\{1\}$. Now $V''K=V$ is tidy
by \cite[Prop.\,3.1]{FLA}.

{\em $\alpha$ is tidy}: Since every compact open subgroup
of~$G$ contains a compact open subgroup satisfying {\bf T1}
by the first step of the algorithm just cited,
and any such is tidy by (d),
we see that~$G$ has small subgroups tidy for~$\alpha$.

(e) Can be proved as \cite[Thm.\,3.32\,(6) and (6$'$)]{BaW}.

(f) To see that $\mu$ is injective,
suppose that $\mu(x,y,z)=\mu(x',y',z')$.
Then $y^{-1}x^{-1}x'y'=z(z')^{-1}\in M_\alpha U_\alpha \cap U_{\alpha^{-1}}=P_\alpha\cap
U_{\alpha^{-1}}=U_0=\{1\}$ by \cite[La.\,3.29]{BaW},
whence $z=z'$ and $xy=x'y'$. Thus $x^{-1}x'=y(y')^{-1}\in U_\alpha\cap M_\alpha=\{1\}$,
using~(c), and hence $x=x'$ and $y=y'$. Thus~$\mu$ is injective.
Let $V\sub G$ be tidy for~$\alpha$.
Then $V=V_-V_+$
where $V_-=(U_\alpha\cap V)V_0$
as a consequence of \cite[Prop.\,3.16]{BaW},
and likewise $V_+=V_0(U_{\alpha^{-1}}\cap V)$.
Here $U_\alpha\cap V$,
$U_{\alpha^{-1}}\cap V$
and $V_0=M_\alpha\cap V$ are compact open subgroups of $U_\alpha$,
$U_{\alpha^{-1}}$ and~$M_\alpha$, respectively.
Since $V=(U_\alpha\cap V_-)V_0(U_{\alpha^{-1}}\cap V)$,
the map~$\mu$ induces a continuous bijection between
$W:=(U_\alpha\cap V)\times
(M_\alpha\cap V)\times (U_{\alpha^{-1}}\cap V)$ and~$V$.
Both~$W$ and~$V$ being compact, we deduce that
$\mu|_W^V$ is a homeomorphism.
To verify that $\im(\mu)$ is open and~$\mu$
a homeomorphism, let $(x_0,y_0,z_0)\in U_\alpha \times M_\alpha \times U_{\alpha^{-1}}$.
As $M_\alpha$ normalizes $U_\alpha$,
the formula $h(x,y,z):=(y_0^{-1}xy_0,y,z)$
defines a homeomorphism $h$ of
$U_\alpha\times M_\alpha\times U_{\alpha^{-1}}$ onto itself.
Left and right translations being homeomorphisms
and $h^{-1}(W)$ being an identity neighbourhood,
we now deduce from
\[
x_0y_0 \mu|_W(h(x,y,z))z_0\; =\; \mu(x_0x, y_0 y, zz_0)\qquad\mbox{for all $\, (x,y,z)\in h^{-1}(W)$}
\]
that $\mu$ takes a neighbourhood of $(x_0,y_0,z_0)$
homeomorphically onto a neighbourhood of~$\mu(x_0,y_0,z_0)$.
\end{proof}
{\footnotesize 
Helge Gl\"{o}ckner, TU Darmstadt, Fachbereich Mathematik AG~5,
Schlossgartenstr.\,7,\\
64289 Darmstadt, Germany. E-Mail: gloeckner\at{}mathematik.tu-darmstadt.de}

\begin{thebibliography}{99}
%
\bibitem{BaW}
Baumgartner, U. and G.\,A. Willis, {\em Contraction groups and scales of automorphisms
of totally disconnected locally compact groups},
to appear in Israel J. Math.\ (cf.\ arXiv:math.GR/0302095).
%
%
\bibitem{DaS}
Dani S.\,G. and R. Shah, {\em Contraction subgroups and semistable
measures on $p$-adic Lie groups}, Math.\ Proc.\ Camb.\
Phil.\ Soc.\ {\bf 110}\,(1991), 299--306.
%
%
\bibitem{SCA}
Gl\"{o}ckner, H., {\em Scale functions on $p$-adic Lie groups},
Manuscr.\ Math.\ {\bf 97}\,(1998), 205--215.
%
%
\bibitem{SC2}
-----, {\em Scale functions on Lie groups over local fields of positive
characteristic}, in preparation (will be posted in the arXiv in September~2004).
%
%
\bibitem{GaW}
Gl\"{o}ckner, H. and G.\,A. Willis, {\em Uniscalar $p$-adic Lie groups},
Forum Math.\ {\bf 13}\,(2001), 413--421.
%
%
\bibitem{Haz}
Hazod, W., {\em Remarks on $($semi-$)$ stable probabilities},
pp.\,182--203 in: Heyer, H. (Ed.),
``Probability Measures on Groups VII,''
Lecture Notes in Math.\ {\bf 1064}, Springer-Verlag, 1984.
%
%
\bibitem{HaS}
Hazod, W. and E. Siebert, {\em Automorphisms on a Lie group contracting modulo a compact subgroup
and applications to semistable convolution semigroups},
J. Theor.\ Probab.\ {\bf 1}\,(1988), 211--225.
%
%
\bibitem{Sie}
Siebert, E., {\em Contractive automorphisms on locally compact
groups}, Math.\ Z. {\bf 191}\,(1986), 73--90.
%
%
\bibitem{Wan}
Wang, J.\,S.\,P., {\em The Mautner phenomenon for $p$-adic
Lie groups}, Math.\ Z. {\bf 185}\,(1984), 403--412.
%
%
\bibitem{Wi1}
Willis, G.\,A., {\em The structure of totally disconnected, locally compact groups},
Math.\ Ann.\ {\bf 300}\,(1994), 341--363.
%
%
\bibitem{Wi2}
-----, {\em Further properties of the scale function on a totally disconnected group},
J. Algebra {\bf 237}\,(2001), 142--164.
%
%
\bibitem{FLA}
-----, {\em Tidy subgroups for commuting automorphisms
of totally disconnected locally compact groups: An analogue of simultaneous
triangularisation of matrices}, New York J. Math.\ {\bf 10}\,(2004),
1--35.
%
%
\end{thebibliography}
\end{document}